\documentstyle[12pt]{article}
\makeatletter \oddsidemargin  -.1in \evensidemargin -.1in
\newcommand{\singlespacing}{\let\CS=\@currsize\renewcommand{\baselinestretch}{1}\tiny\CS}
\newcommand{\oneandahalfspacing}{\let\CS=\@currsize\renewcommand{\baselinestretch}{1.25}\tiny\CS}
\newcommand{\doublespacing}{\let\CS=\@currsize\renewcommand{\baselinestretch}{1.35}\tiny\CS}

\newtheorem{rule-def}[theorem]{Rule}

\RequirePackage[dvips]{graphicx} \textheight 23.0 cm

\begin{document}

\title{\bf Mathematical modelling of blood flow through a tapered overlapping stenosed artery with variable viscosity}
\author{\small G. C. Shit\thanks{Correspondence author.
{\it Email address:} gcs@math.jdvu.ac.in (G. C. Shit) Tel: +91 3324572716}, ~{\small M. Roy}
~and {\small~A. Sinha}\\
\it \small Department of Mathematics, Jadavpur University,
Kolkata-700032, India \\ }
\date{}
\maketitle
\doublespacing
\begin{abstract}
This paper presents a theoretical study of blood flow through a tapered and overlapping stenosed artery under the action of an externally applied magnetic field. The fluid (blood) medium is assumed to be porous in nature. The variable viscosity of blood depending on hematocrit (percentage volume of erythrocytes) is taken into account in order to improve resemblance to the real situation. The governing equation for laminar, incompressible and Newtonian fluid subject to the boundary conditions is solved by using a well known Frobenius method. The analytical expressions for velocity component, volumetric flow rate, wall shear stress and pressure gradient are obtained. The numerical values are extracted from these analytical expressions and are presented graphically. It is observed that the influence of hematocrit, magnetic field and the shape of artery have important impact on the velocity profile, pressure gradient and wall shear stress. Moreover, the effect of primary stenosis on the secondary one has been significantly observed.

\noindent {\bf Keywords:} {\it Overlapping Stenoses, Tapered Artery, MHD Flow, Porous Vessel, Hematocrit, Frobenius Method}

\end{abstract}
\section {Introduction}
Many cardiovascular diseases such as due to the arterial occlusion is one of the leading cause of death world wide. The partial occlusion of the arteries due to stenotic obstruction not only restrict the regular blood flow but also characterizes the hardening and thickening of the arterial wall. However, the main cause of the formation of stenosis is still unknown but it is well established that the fluid dynamical factors play an important role as to further development of stenosis. Therefore, during the past few decay several studies were conducted by  Young \cite{R1}, Young and Tsai \cite{R2, R3} to understand the effects of stenosis on blood flow through arteries.
\begin{center}
\begin{tabular}{|l l|}\hline
 {~\bf Nomenclature} &~\\
  $z$ & Axial distance\\
  $r$ & Radial distance\\
  $R(z)$ & Arterial wall\\
  $d$ & Distance from inlet of onset of stenosis\\
  $u$ & Axial velocity component\\
  $\bar{u}$ & Non-dimensional axial velocity\\
  $\mu(r)$ & Coefficient of viscosity of blood at a radial distance $r$\\
  $\mu_0$ & Viscosity coefficient for plasma\\
  $h(r)$ & Hematocrit at a distance $r$\\
  $H$ & Maximum hematocrit at the center of the arterial segment\\
  $m({\ge} 2)$ & Shape parameter of hematocrit\\
  $R_0$ & Radius of the pericardial surface of the normal portion of the arterial segment\\
  $l$ & Length between throat of two stenosis stenoses\\
  $L$ & Length of arterial segment\\
  $B_0$ & Magnetic field strength\\
  $\sigma$ & Electrical conductivity\\
  $k$ & Porous permeability parameter\\
  $M$ & Hartmann number\\
  $p$ & Blood pressure\\
  $\bar{\frac{dp}{dz}}$ & Non-dimensional axial pressure gradient\\
  $Q$ & Volumetric flow rate\\
  $\bar{Q}$ & Non-dimensional volumetric flow rate\\
  $\tau_R$ & Wall Shear stress\\
  $\bar{\tau}$ & Non-dimensional wall shear stress\\
   \hline
\end{tabular}
\end{center}

Tu and Deville \cite{R4} investigated pulsatile flow of blood in stenosed arteries. Misra and Shit \cite{R5, R6} studied in two different situations on the blood flow through arterial stenosis by treating blood as a non-Newtonian ( Herschel-Bulkley ) fluid model. It is generally well known that blood, being a suspension of red cells in plasma, behaves like a non-Newtonian fluid at low shear rates \cite{R7,R8}. However, Misra and Chakravarty \cite{R9} and Shit and Roy \cite{R10} put forwarded a mathematical analysis for the unsteady flow of blood through arteries having stenosis, in which blood was treated as a Newtonian viscous incompressible fluid.

The hemodynamics associated with a single stenotic lesion are significantly affected by the presence of a second lesion. In many situations there are evidences of the occurrence of the multiple or overlapping stenosis such as the patients of angiograms. Misra et al. \cite{R11} conducted a theoretical study for the effects of multiple stenosis. An experimental study of blood flow through an arterial segment having multiple stenoses were made by Talukder et al. \cite{R12}. The effects of overlapping stenosis through an arterial stenosis have been successfully carried out analytically as well as numerically by Chakravarty and Mandal \cite{R13} and Layek et al. \cite{R14} respectively. However, all these studies are restricted in the consideration of magnetic field and the porous medium.

Since blood is electrically conducting fluid, its flow characteristics is influenced by the application of magnetic field. If a magnetic field is applied to a moving and electrically conducting fluid, it will induce electric as well as magnetic fields. The interaction of these fields produces a body force per unit volume known as Lorentz force, which has significant impact on the flow characteristics of blood. Such an analysis may be useful for the reduction of blood flow during surgery and Magnetic Resonance Imaging ( MRI ). The effect of magnetic field on blood flow has been analyzed theoretically and experimentally by many investigators \cite{R15} - \cite{R18} under different situations. Shit and his co-investigators \cite{R19} - \cite{R23} explored variety of flow behaviour of blood in arteries by treating Newtonian/ non-Newtonian model in the presence of a uniform magnetic field.

Very recently, the studies of blood flow through porous medium has gained considerable attention to the medical practitioners/ clinicians because of its enormous changes in the flow characteristics. The capillary endothelium is, in turn, covered by a thin layer lining the alveoli, which has been treated as a porous medium. Dash et al. \cite{R24} considered the Brinkman equation to model the blood flow when there is an accumulation of fatty plaques in the lumen of an arterial segment and artery-clogging takes place by blood clots. They considered the clogged region as a porous medium. Bhargava et al. \cite{R25}  studied the transport of pharmaceutical species in laminar, homogeneous, incompressible, magneto-hydrodynamic, pulsating flow through two-dimensional channel with porous walls containing non-Darcian porous materials. Misra et al. \cite{R20, R26} presented a mathematical model as well as numerical model for studying blood flow through a porous vessel under the action of magnetic field, in which the viscosity varies in the radial direction.

Hematocrit is the most important determinant of whole blood viscosity. Therefore, blood viscosity and vascular resistance ( due to the presence of stenosis ) affect total peripheral resistance to blood flow, which is abnormally high in the primary stage of hypertension. Again hematocrit is a blood test that measures the percentage of red blood cells present in the whole blood of the body. The percentage of red cells in adult human body is approximately 40 - 45 \% \cite{R27}. Red cells may affect the viscosity of whole blood and thus the velocity distribution depends on the hematocrit. So blood can not be considered as homogeneous fluid \cite{R14}. Due to the high shear rate near the arterial wall, the viscosity of blood is low and the concentration of cells is high in the central core region. Therefore, blood may be treated as Newtonian fluid with variable viscosity particularly in the case of large blood vessels.

The present study is motivated towards a theoretical investigation of blood flow through a tapered and overlapping stenosed artery in the presence of magnetic field. The study pertains to a situation in which the variable viscosity of blood depending upon hematocrit is taken into consideration. The present model is designed in such a way that it could be applicable to both converging/ diverging artery depending on the choice of tapering angle $\alpha$. Thus, the study will answers the question of mechanism of further deposition of plaque under various aspects.

\section {Mathematical modelling of the problem }
We consider the laminar, incompressible and Newtonian flow of blood through axisymmetric two-dimensional tapered and overlapping stenosed artery. Any material point in the fluid is representing by the cylindrical polar coordinate $(r, \theta, z)$, where $z$  measures along the axis of the artery and that of $r$ and $\theta$ measure along the radial and circumferential directions respectively. The mathematical expression that corresponds to the geometry of the present problem is given by
\begin{flushleft}
\begin{eqnarray}
  R(z)&=&R_{0}\Big[1.0-\frac{11.0}{32.0}l^{3}(z-d)+\frac{47.0}{48.0}l^{2}(z-d)^{2}-l(z-d)^{3}+(1.0/3.0)(z-d)^{4}\Big]a(z)
  \nonumber\\
    &&\hspace{3.5 cm} where~~ d\leq z \leq d+\frac{3 l}{2}\nonumber\\
        &=&R_{0}a(z)  ~~~~~~~~~~~~~~~~~~~~~ elsewhere
\end{eqnarray}
\end{flushleft}
where the onset of the stenosis is located at a distance $d$ from the inlet, $\frac{3 l}{2}$ the length of the overlapping stenosis and $l$ representing the distance between two critical height of the stenoses. The expression for $a(z)$ is responsible for the artery to be converging or diverging depending on tapering angle $\alpha$ has the form
\begin{eqnarray}
a(z) = 1+z~tan(\alpha).
\end{eqnarray}
\begin{figure}[htbp]
\begin{center}
        \includegraphics[width=4.8in,height=3.0in ] {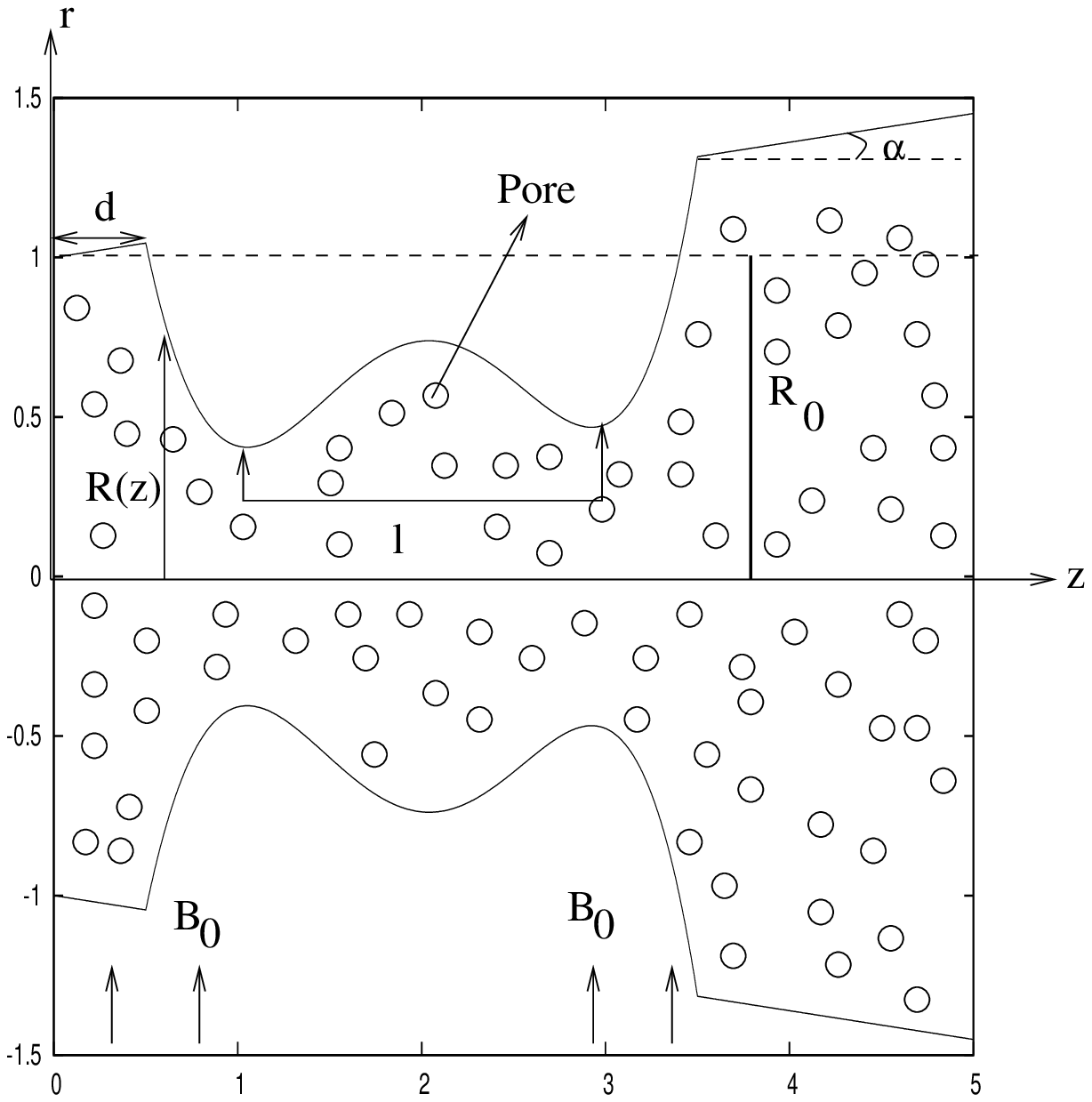}\\
                   Fig. 1 Schematic diagram of the model geometry.\\
\end{center}
\end{figure}\\
We assumed that blood is incompressible, suspension of erythrocytes in plasma and has uniform dense throughout but the viscosity $\mu(r)$ varies in the radial direction. According to Einstein's formula for the variable viscosity of blood taken to be
\begin{eqnarray}
 \mu(r)=\mu_{0}\left[1+\beta h(r)\right],
\end{eqnarray}
where $\mu_0$ is the coefficient of viscosity of plasma, $\beta$ is a constant (whose value for blood is equal to 2.5) and $h(r)$ stands for the hematocrit. The analysis will be carried out by using the following empirical formula for hematocrit given by in Lih \cite{R29}
\begin{equation}
h(r)=H\Big[1-\Big(\frac{r}{R_0}\Big)^m \Big],
\end{equation}\\
in which $R_0$ represents the radius of a normal arterial segment, $H$ is the maximum hematocrit at the center of the artery and $m(\ge 2)$ a parameter that determines the exact shape of the velocity profile for blood. The shape of the profile given by Eq.(4) is valid only for very dilute suspensions of erythrocytes, which are considered to be of spherical shape.

According to our considerations, the equation that governs the flow of blood under the action of an external magnetic field through porous medium may be put as
\begin{equation}
\frac{\partial p}{\partial z}-\frac{1}{r} \frac{\partial(r \mu(r)\frac{du}{d r})}{\partial r}
+\sigma B_0^2u+\frac{\mu(r)}{\bar{k}}u=0,
\end{equation}\\
where $u$ denotes the (axial) velocity component of blood, $p$ the blood pressure, $\sigma$ the electrical conductivity, $\bar{k}$ the permeability of the porous medium and $B_0$ is the applied magnetic field strength.

To solve our problem, we use no-slip boundary condition at the arterial wall, that is,
\begin{eqnarray}
 u=0 ~~~~ at ~~~~~r = R(z).
\end{eqnarray}
Further we consider axi-symmetric boundary condition of axial velocity at the mid line of the artery as
\begin{eqnarray}
 \frac{\partial u}{\partial r}=0 ~~~~ at ~~~~~r = 0
\end{eqnarray}

\section {Analytical solution}
In order to simplify our problem, let us introduce the following transformation
\begin{equation}
\xi=\frac{r}{R_0}
\end{equation}
With the use of the transform defined in (8) and the equations (3) and (4), the governing equation (5) reduces to
\begin{equation}
\frac{1}{\xi}\frac{\partial}{\partial \xi}[\xi(a_1-a_2\xi^m)\frac{\partial u}{\partial \xi}]-M^2u-\frac{1}{k}(a_1-a_2\xi^m)u=\frac{R_0^2}{\mu_0}\frac{\partial p}{\partial z}
\end{equation}
with~~$a_1=1+a_2$,~ $a_2=\beta H$,~~$k=\frac{\bar{k}}{R_0^2}$~ and~ $M^2=\frac{\sigma (B_0R_0)^2}{\mu_0}$.\\
Similarly the boundary conditions transformed into
\begin{equation}
u=0~~~~~~~~~~~~~~~~~~~~~~~~~~~  at~~~ \xi=\frac{R(z)}{R_0}
\end{equation}
\begin{equation}
and~~~~\frac{\partial u}{\partial \xi}=0 ~~~~~~~~~~~~~~~~~~~~~~~~~~~~~~at ~~~\xi=0
\end{equation}
The equation (9) can be solved subjected to the boundary conditions (10) and (11) using Frobenius method. For this, of course, $u$ has to be bounded at $\xi = 0$. Then only admissible series solution of the equation (9) will exists and can put in the form
\begin{equation}
u=K\sum\limits_{i=0}^{\infty}
{A_i\xi^i}+\frac{R_0^2\frac{dp}{dz}}{4a_1\mu_0}\sum\limits_{i=0}^{\infty}
{B_i\xi^{i+2}}~,
\end{equation}
where $K$, $A_i$ and $B_i$ are arbitrary constants.\\
To find the arbitrary constant $K$, we use the no-slip boundary condition (10) and obtained as
 \begin{equation}
K=-\frac{R_0^2 \frac{dp}{dz}}{4a_1\mu_0}
\frac{\sum\limits_{i=0}^{\infty}
{B_i(\frac{R}{R_0})^{i+2}}}{\sum\limits_{i=0}^{\infty}{A_i(\frac{R}{R_0})^{i}}}.
\end{equation}
Substituting the value of $u$ from equation (12) into equation (9) and we get
\begin{eqnarray}
K\Big[\sum_{i=0}^{\infty}i(i-1)(a_{1}-a_{2}\xi^{m})A_{i}\xi^{i-2}+\sum_{i=0}^{\infty}i(a_{1}-(m-1)a_{2}\xi^{m})A_{i}\xi^{i-2}-(M^{2}+\frac{a_{1}}{k})\sum_{i=0}^{\infty}A_{i}\xi^{i}\nonumber\\ +\sum_{i=0}^{\infty}\frac{a_{2}}{k}A_{i}\xi^{i+m}\Big]+\frac{R_{0}^{2}\frac{dp}{dz}}{4a_{1}\mu_{0}}\Big[\sum_{i=0}^{\infty}(i+1)(i+2)(a_{1}-a_{2}\xi^{m})B_{i}\xi^{i}+\sum_{i=0}^{\infty}(i+2)(a_{1}-(m-1)a_{2}\xi^{m})B_{i}\xi^{i}\nonumber\\
-(M^{2}+\frac{a_{1}}{k})\sum_{i=0}^{\infty}B_{i}\xi^{i+2}+\sum_{i=0}^{\infty}\frac{a_{2}}{k}B_{i}\xi^{i+m+2}\Big]=\frac{R_{0}^{2}\frac{dp}{dz}}{4a_{1}\mu_{0}}
\end{eqnarray}
Equating the coefficients of $K$ and other part in equation (14) we have,
\begin{eqnarray}
 \sum_{i=0}^{\infty}i(i-1)(a_{1}-a_{2}\xi^{m})A_{i}\xi^{i-2}+\sum_{i=0}^{\infty}i(a_{1}-(m-1)a_{2}\xi^{m})A_{i}\xi^{i-2}-(M^{2}+\frac{a_{1}}{k})\sum_{i=0}^{\infty}A_{i}\xi^{i} \nonumber\\ +\sum_{i=0}^{\infty}\frac{a_{2}}{k}A_{i}\xi^{i+m}=0
\end{eqnarray}
and
\begin{eqnarray}
\sum_{i=0}^{\infty}(i+1)(i+2)(a_{1}-a_{2}\xi^{m})B_{i}\xi^{i}+\sum_{i=0}^{\infty}(i+2)(a_{1}-(m-1)a_{2}\xi^{m})B_{i}\xi^{i}\nonumber\\-(M^{2}+\frac{a_{1}}{k})\sum_{i=0}^{\infty}B_{i}\xi^{i+2}+\sum_{i=0}^{\infty}\frac{a_{2}}{k}B_{i}\xi^{i+m+2}=1
\end{eqnarray}
Hence the constants $A_i$ and $B_i$ are obtained by equating the coefficients of $\xi^{i+1}$ and $\xi^{i-1}$ from both side of equations (15) and (16) respectively and can be put in the form
\begin{equation}
A_{i+1}=\frac{a_2(1+i)^{2}A_{i+1-m}-ma_2A_{i+1-m}+(M^2+\frac{a_1}{k})A_{i-1}-\frac{a_2}{k}A_{i-1-m}}{a_1(1+i)^2}
\end{equation}
\begin{equation}
B_{i+1}=\frac{a_2(3+i)^{2}-a_{2}m B_{i+1-m}+(M^2+\frac{a_1}{k})B_{i-1}-\frac{a_2}{k}B_{i-1-m}}{a_1(3+i)^2}
\end{equation}
with
\begin{equation}
A_0=B_0=1.
\end{equation}
 Substituting the expression for $K$ in the equation (12), we have
\begin {equation}
u=-\frac{R_0^2 \frac{dp}{dz}}{4a_1\mu_0} \frac{\Bigg
[\sum\limits_{i=0}^{\infty}
{B_i(\frac{R}{R_0})^{i+2}}\sum\limits_{i=0}^{\infty}{A_i
{\xi}^{i}}-\sum \limits_{i=0}^{\infty} {B_i {\xi}^{i+2}}\sum
\limits_{i=0}^{\infty}{A_i(\frac{R}{R_0})^{i}}\Bigg]}{\sum\limits_{i=0}^{\infty}{A_i(\frac{R}{R_0})^{i}}}
\end {equation}
The average velocity $u_0$ has the form
\begin {equation}
u_0=-\frac{R_0^2}{8\mu_0}\Bigg(\frac{dp}{dz}\Bigg)_0
\end {equation}
where $\Bigg(\frac{dp}{dz}\Bigg)_0$ is the pressure gradient of the flow field in the normal artery in the absence of magnetic field. The non-dimensional expression for $u$ is given by
\begin {equation}
\bar u=\frac{u}{u_0}=\frac{2}{a_1}\frac{\frac{dp}{dz}}{(\frac{dp}{dz})_0}
\frac{\Bigg [\sum\limits_{i=0}^{\infty}
{B_i(\frac{R}{R_0})^{i+2}}\sum\limits_{i=0}^{\infty}{A_i
{\xi}^{i}}-\sum \limits_{i=0}^{\infty} {B_i {\xi}^{i+2}}\sum
\limits_{i=0}^{\infty}{A_i(\frac{R}{R_0})^{i}}\Bigg]}{\sum\limits_{i=0}^{\infty}{A_i(\frac{R}{R_0})^{i}}}
\end {equation}
The volumetric flow rate across the arterial segment is given by
\begin{equation}
Q=\int \limits_{0}^{\frac{R}{R_0}}2\pi R_0\xi u(\xi)\,d\xi
\end{equation}
Substituting $u$ from Eq.(20) into Eq.(23) and then integrating with respect to $\xi$, we obtain
\begin{equation}
Q=-\frac{\pi R_0^3 \frac{dp}{dz}}{2a_1\mu_0}\frac{\Bigg
[\sum\limits_{i=0}^{\infty}
{B_i(\frac{R}{R_0})^{i+2}}\sum\limits_{i=0}^{\infty}\frac{{A_i
{(\frac{R}{R_0})}^{i+2}}}{(i+2)}- \sum \limits_{i=0}^{\infty}
\frac {{B_i {(\frac{R}{R_0})}^{i+4}}}{(i+4)}\sum
\limits_{i=0}^{\infty}{A_i(\frac{R}{R_0})^{i}}\Bigg]}{\sum\limits_{i=0}^{\infty}{A_i(\frac{R}{R_0})^{i}}}
\end{equation}
If $Q_0$ be the volumetric flow rate in the normal portion of the artery, in the absence of magnetic field and porosity effect, then
\begin{equation}
Q_0=-\frac{\pi R_0^3}{8\mu_0}\big(\frac{dp}{dz}\big)_0
\end{equation}
Therefore, the non-dimensional volumetric flow rate has the following form
\begin{equation}
\bar Q=\frac{Q}{Q_0}=\frac{4 \frac{dp}{dz}}{a_1 (\frac{dp}{dz})_0}\frac{\Bigg
[\sum\limits_{i=0}^{\infty}
{B_i(\frac{R}{R_0})^{i+2}}\sum\limits_{i=0}^{\infty}\frac{{A_i
{(\frac{R}{R_0})}^{i+2}}}{(i+2)}- \sum \limits_{i=0}^{\infty}
\frac {{B_i {(\frac{R}{R_0})}^{i+4}}}{(i+4)}\sum
\limits_{i=0}^{\infty}{A_i(\frac{R}{R_0})^{i}}\Bigg]}{\sum\limits_{i=0}^{\infty}{A_i(\frac{R}{R_0})^{i}}}
\end{equation}
If the flow is steady and no outward/inward flow takes place through the arterial segment, then the mass flux is constant and hence $\bar Q=1$. The expression for pressure gradient from (26) can be put as
\begin{equation}
(\frac{\overline{dp}}{dz})=\frac{\frac{dp}{dz}}{(\frac{dp}{dz})_0}=\frac{a_1}{4}\frac{\sum
\limits_{i=0}^{\infty} {A_i(\frac{R}{R_0})^{i}}}{\Bigg
[\sum\limits_{i=0}^{\infty}
{B_i(\frac{R}{R_0})^{i+2}}\sum\limits_{i=0}^{\infty}\frac{{A_i
{(\frac{\eta}{R_0})}^{i+2}}}{(i+2)}- \sum \limits_{i=0}^{\infty}
\frac {{B_i {(\frac{R}{R_0})}^{i+4}}}{(i+4)}\sum
\limits_{i=0}^{\infty}{A_i(\frac{R}{R_0})^{i}}\Bigg]}
\end{equation}
The wall shear stress on the endothelial surface is given by
\begin{equation}
\tau_{R}=\Bigg[-\mu(r)\frac{du}{d r}\Bigg]_{r=R(z)}
\end{equation}
Substituting $u$ from Eq.(20) into Eq.(28), we obtain
\begin{equation}
\tau_{R}=\frac{\frac{dp}{dz}R_0 \Big[1+\beta H\{1-\Big(\frac{R}{R_0}\Big)^m\}\Big]}{4a_1}\frac{\Bigg [{\sum\limits_{i=0}^{\infty}
{B_i(\frac{R}{R_0})^{i+2}}\sum\limits_{i=0}^{\infty}{iA_i
{(\frac{R}{R_0})}^{i-1}}}- \sum \limits_{i=0}^{\infty} {(i+2)B_i
{(\frac{R}{R_0})}^{i+1}}\sum
\limits_{i=0}^{\infty}{A_i(\frac{R}{R_0})^{i}}\Bigg]}{\sum\limits_{i=0}^{\infty}{A_i(\frac{R}{R_0})^{i}}}
\end{equation}
If $\tau_N=-\frac{R_0}{2}(\frac{dp}{dz})_0$ be the shear stress at the normal portion of the arterial
wall, in the absence of magnetic field, the non-dimensional form of the wall shear stress is given by
\begin{eqnarray}
\bar \tau&=&\frac{\tau_{R}}{\tau_N}\nonumber \\
&=&\frac{1}{2a_1}(\frac{\overline {dp}}{dz})[1+\beta
H(1-(\frac{R}{R_0})^m)]\frac{\frac{dp}{dz}}{\sum\limits_{i=0}^{\infty}{A_i(\frac{R}{R_0})^{i}}}
\Bigg [{\sum\limits_{i=0}^{\infty}
{B_i(\frac{R}{R_0})^{i+2}}\sum\limits_{i=0}^{\infty}{iA_i {(\frac{R}{R_0})}^{i-1}}}\nonumber\\
& &\hspace{1.0 cm}-\sum \limits_{i=0}^{\infty} {(i+2)B_i
{(\frac{R}{R_0})}^{i+1}}\sum
\limits_{i=0}^{\infty}{A_i(\frac{R}{R_0})^{i}}\Bigg].
\end{eqnarray}~~~~~~~~

\section{Results and Discussion}
In the previous section we have obtained analytical expressions for different flow characteristics of blood through porous medium under the action of an external magnetic field. In this section we are to discuss the flow characteristics graphically with the use of following valid numerical data which is applicable to blood. To continue this section we have use the following standard values of physical parameter:
\[ l=2.0,~ d=0.5,~ \alpha = 0.09,~ H=0.2,~ k=0.25,~ \beta= 2.5,~ L=5,~ M=2.5,\]
Fig. 2 indicates the different locations of the stenosis in the axial direction. In Fig. 2, $z=0.5$ and $z=3.5$ correspond to the onset and outset of the stenosis and $z=1$ and $z=3$ represent the throat of the primary stenosis and secondary stenosis respectively. It is interesting to note from this figure that $z=2$ is the location where further deposition takes place and hence it is known as overlapping stenosis.

\begin{figure}[htbp]
\begin{center}
        \includegraphics[width=4.8in,height=3.0in ] {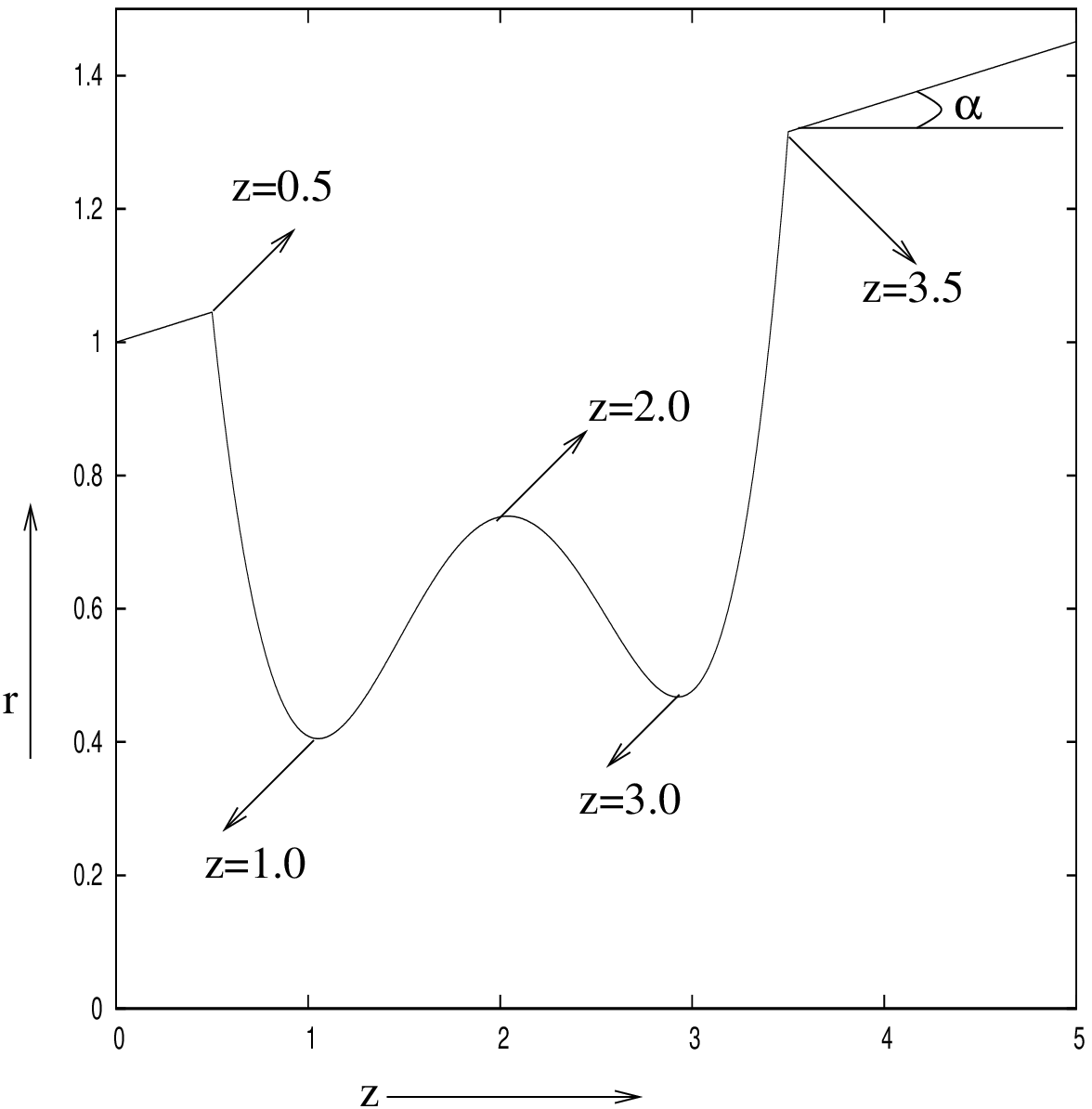}\\
                   Fig. 2 Locations of the points z=0.5, 1.0, 2.0, 3.0, 3.5.\\
\end{center}
\end{figure}
The variation of axial velocity at different axial position along the radial direction are shown in Fig. 3. We observe from this figure that the velocity is maximum at the central line of the vessel for all position of $z$. Among all these positions, the velocity is high at the throat of the primary stenosis and low at the onset of the overlapping stenosis. However, the central line velocity at the throat of secondary stenosis is about 30 \% less than the primary stenosis. But the central line velocity at $z=2$ (between the throat of two stenosis) suddenly falls about 55 \% than that of secondary one. This observation may leads to the flow circulation zone and may causes further deposition of plaque. Fig. 4 depicts the variation of axial velocity for different length between the throat of two stenoses. It has been observed that the velocity at the throat of the stenoses significantly increases with the increase of $l$. Thus, the effect of the shape of stenosis has important role on the flow characteristics. Similar is the observation from Fig. 5 that the axial velocity decreases as the tapering angle $\alpha$ increases. Fig. 6 illustrates the variation of axial velocity at the throat of the secondary stenosis for different values of the Hartmann Number $M$. We observe that the axial velocity significantly decreases at the central line of the artery with the increase of the magnetic field strength. While the velocity in the vicinity of the arterial wall increases with the increasing values of $M$ in order to maintain constant volumetric flow rate. It is also well known that when a magnetic field is applied in an electrically conducting fluid (here for blood) there arises Lorentz force, which has a tendency to slow down the motion of the fluid. It has been observed from Fig. 7 that the axial velocity near the central line of the channel increases with the increase of the permeability parameter $k$, while the trend is reversed in the vicinity of the arterial wall. This phenomena is noticed because of the permeability parameter $k$ is depend as the reciprocal of the permeability of the porous medium. Fig. 8 gives the distribution of axial velocity for different values of the hematocrit $H$. We note from Fig. 8 that the axial velocity decreases at the core region of the artery with the increase of hematocrit level $H$, whereas the opposite trend is observed in the peripheral region. This fact lies with in the hematocrit as the blood viscosity is high in the core region due to the aggregation of blood cells rather than low viscosity in the plasma near the arterial wall.

Figs. 9 - 12 illustrate the variation of pressure gradient $\bar{\frac{d p}{d z}}$ along the length of the stenosis for different values of the physical parameters of interest. Fig. 9 shows that the axial pressure gradient increases with the increase of magnetic field strength. We have already observed that the Lorentz force has reducing effect of blood velocity, so as more pressure is needed to pass the same amount of fluid under the action of an external magnetic field.  However, the opposite trend is observed in the case of porous permeability parameter $k$ as shown in Fig. 10. It has been seen from Fig. 11 that the pressure gradient increases with the increase of the hematocrit $H$. It indicates from this figure that when the aggregation of blood cells increase at the core region that is hematocrit $H$ is high, more pressure gradient is needed to pass the same amount of the fluid through the stenotic region. It is interesting to note from these three figures that the magnitude of the pressure gradient is high enough at the throat of the primary stenosis than that of the secondary one. But from Fig. 12, we observed that at the throat of secondary stenosis, the pressure gradient is high in comparison to the throat of the primary stenosis. This happens due to the increasing of the tapering angle $\alpha$. Therefore, in the case of diverging artery more pressure is needed as the flow advances in the downstream direction.

Figs. 13 and 14 give the distribution of the wall shear stress for different values of the hematocrit $H$ and tapering angle $\alpha$. We observe from Fig. 13 that the wall shear stress increases as the hematocrit $H$ increases. One can note from this figure that the wall shear stress is low at the throat of the secondary stenosis as well as  at the downstream of the artery. It is generally well known that at the low shear stress region mass transportation takes place and thereby occurs further deposition. However, it is interesting to note from Fig. 14 that the wall shear stress decreases significantly with the increasing values of the tapering angle $\alpha$. Moreover, the magnitude of the wall shear stress is same at both the throat of the stenosis in the absence of tapering angle $\alpha$. Therefore, we may conclude that there is a chance of further deposition at the downstream of the diverging artery.

\section{Conclusions}
A theoretical study of blood flow through overlapping stenosis in the presence of magnetic field has been carried out. In this study the variable viscosity of blood depending on hematocrit and the has been treated as the porous medium. The problem is solved analytically by using Frobenius method. The effects of various key parameters including the tapering angle $\alpha$, percentage of hematocrit $H$, the magnetic field and permeability parameter $k$ are examined. The main findings of the present study may be listed as follows:\\
\begin{itemize}
  \item  The effect of primary stenosis on the secondary one is significant in case of diverging artery $(\alpha >0)$
  \item  The flow velocity at the central region decreases gradually with the increase of magnetic field strength.
  \item  The permeability parameter $k$ has an enhancing effect on the flow characteristics of blood.
  \item  At the core region, the axial velocity decreases with the increase of the percentage of hematocrit $H$.
  \item  The hematocrit and the blood pressure has a linear relationship as reported in \cite{R28}.
  \item  The lower range of hematocrit may leads to the further deposition of cholesterol at the endothelium of the vascular wall.
\end{itemize}
 Finally we can conclude that further potential improvement of the model are anticipated. Since the hematocrit positively affects blood pressure, further study should examine the other factors such as diet, tobacco, smoking, overweight etc. from a cardiovascular point of view. More over on the basis of the present results, it can be concluded that the flow of blood and pressure can be controlled by the application of an external magnetic field.

{\bf Acknowledgement:} {The authors are thankful to Jadavpur University for the financial support through JU Research grant during this investigation.}

\newpage
\vspace{45pt}
\begin{figure} [!ht]
\begin{center}
             \includegraphics[width=4.6in,height=3.6in ] {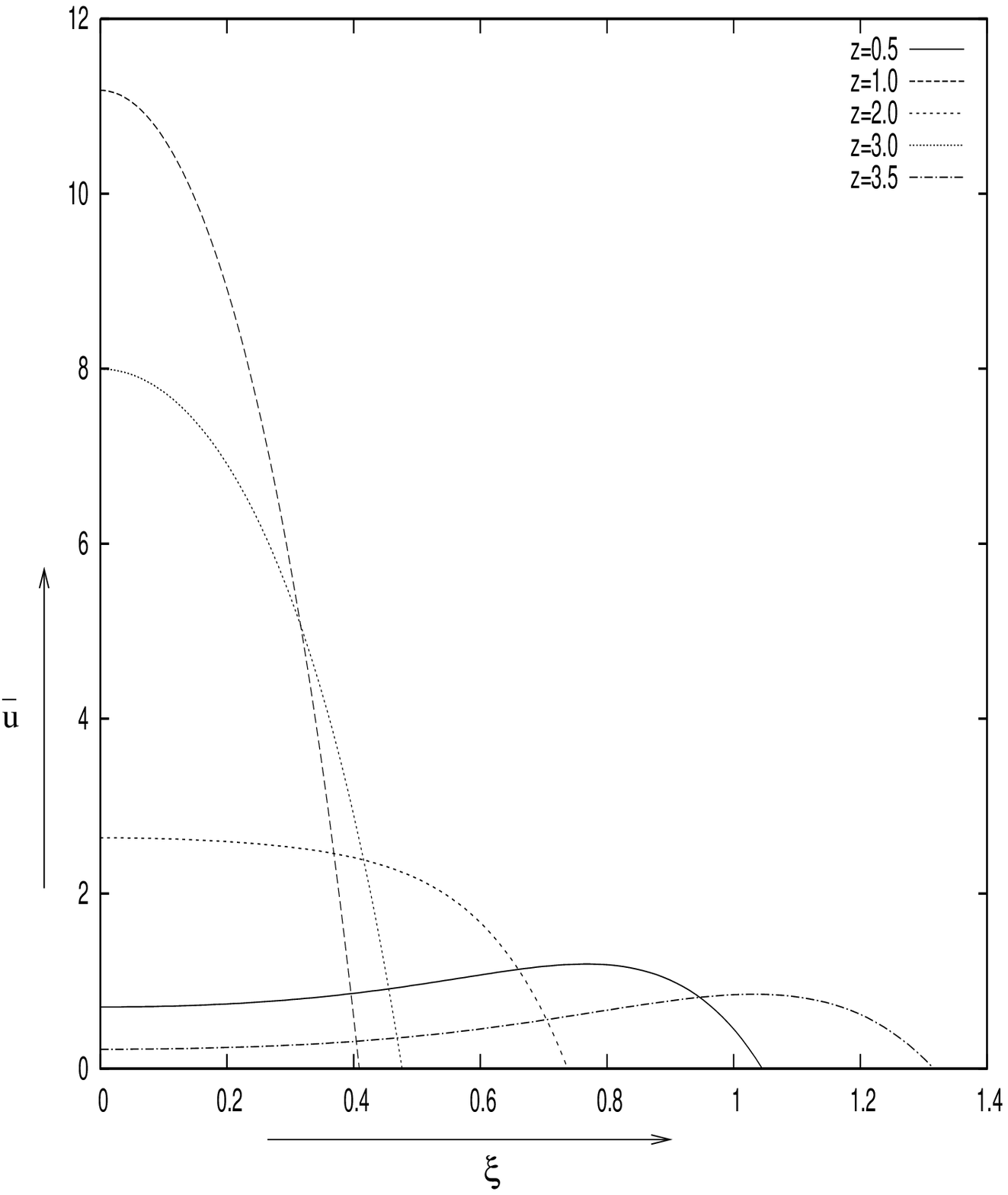}\\
      Fig. 3 Velocity distribution in the radial direction at different axial position $z$, when $H=0.2$, $M=2.5$, $\beta=2.5$, $\alpha=0.09$ and $k=0.25$.\\
            \includegraphics[width=4.6in,height=3.5in ]{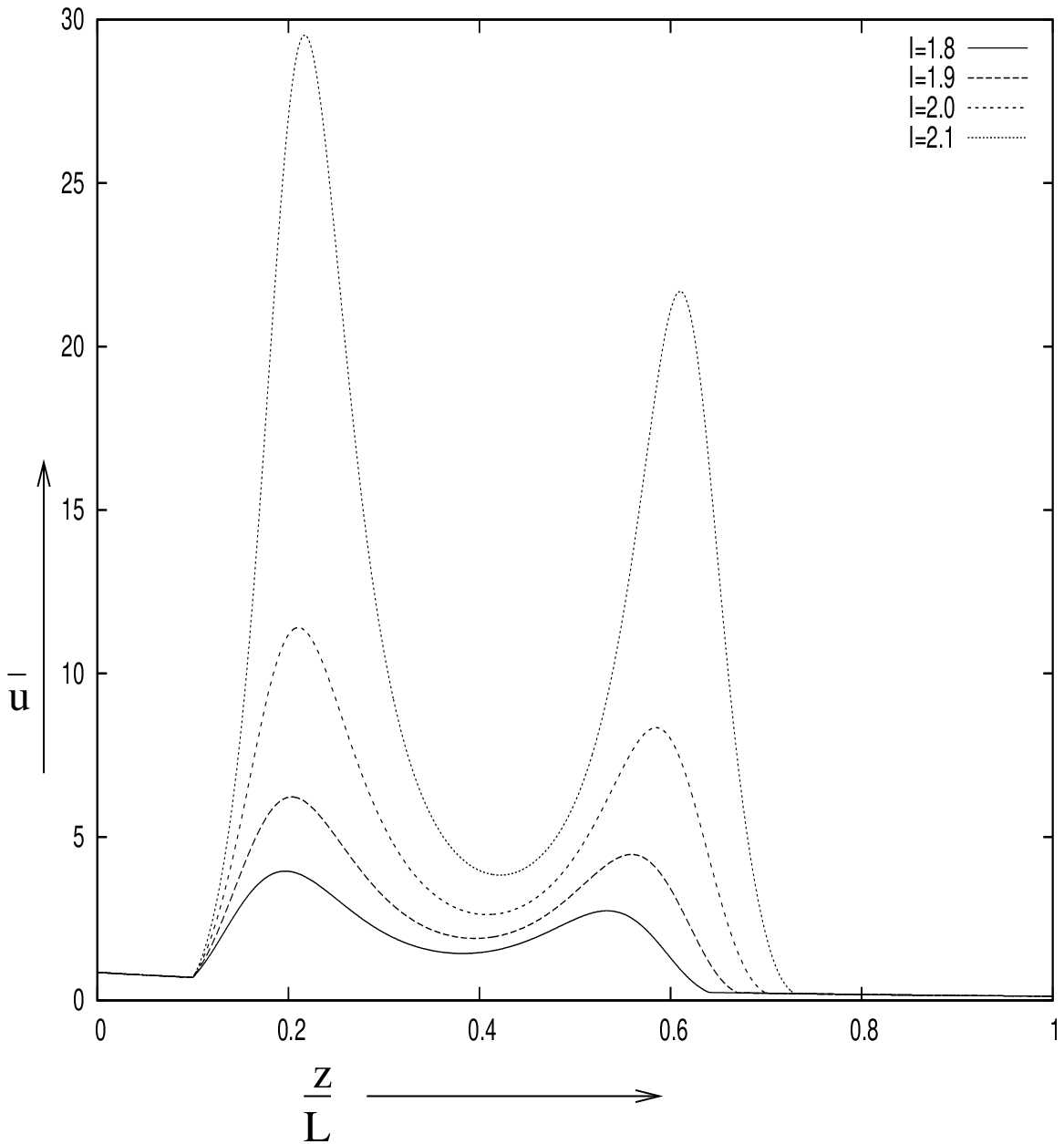}\\
      Fig. 4  Variation of axial velocity along axial direction for different values of $l$, when $H=0.2$, $M=2.5$, $\alpha=0.09$ and $k=0.25$.\\
\end{center}
\end{figure}
\begin{figure}[!ht]
\begin{center}
             \includegraphics[width=4.6in,height=3.5in ]{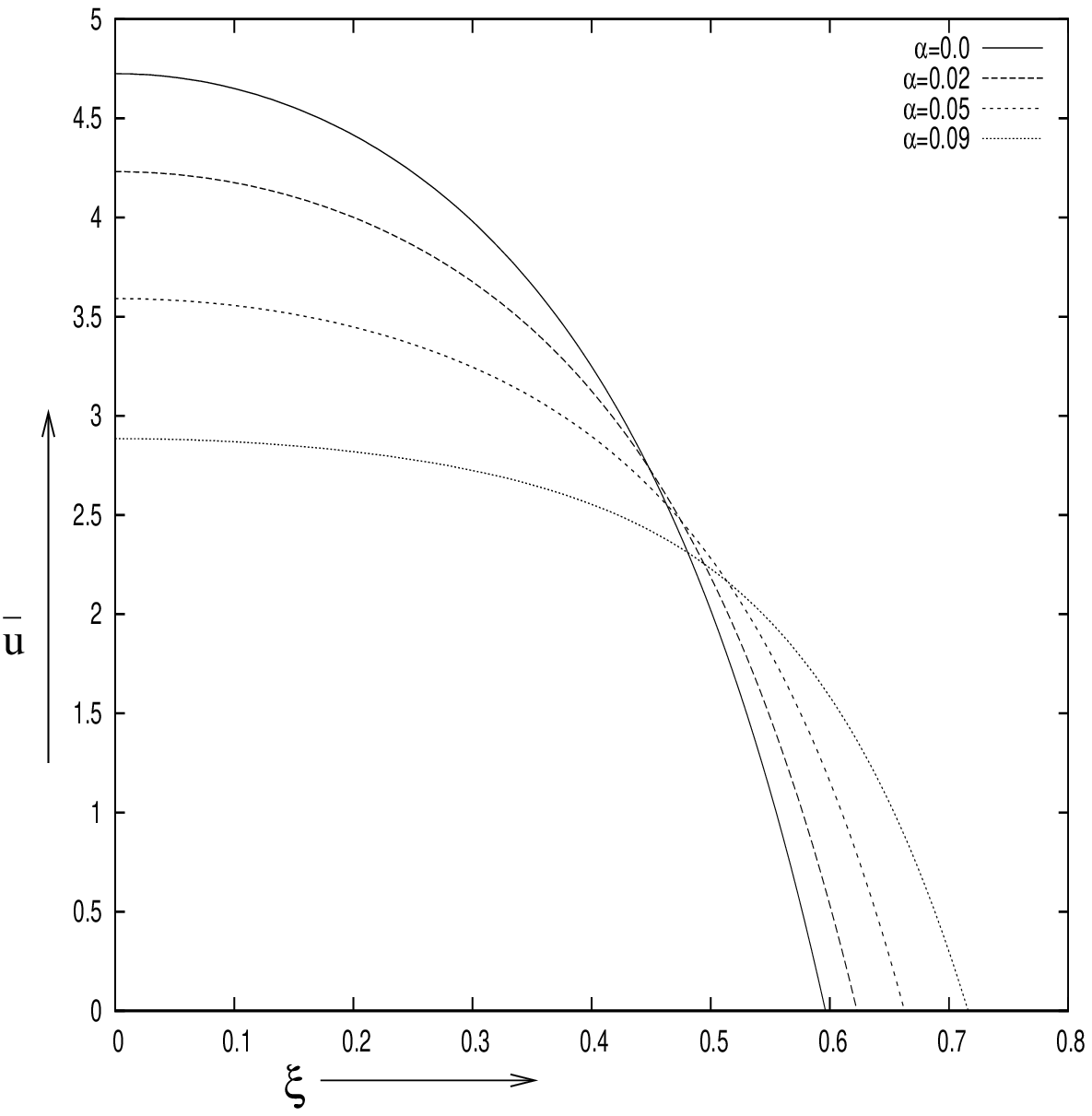}\\
      Fig. 5 Variation of axial velocity  at $z=2.0$ for different values of $\alpha$, when $H=0.2$, $M=2.5$, $\beta=2.5$, and $k=0.25$.\\
             \includegraphics[width=4.6in,height=3.5in ]{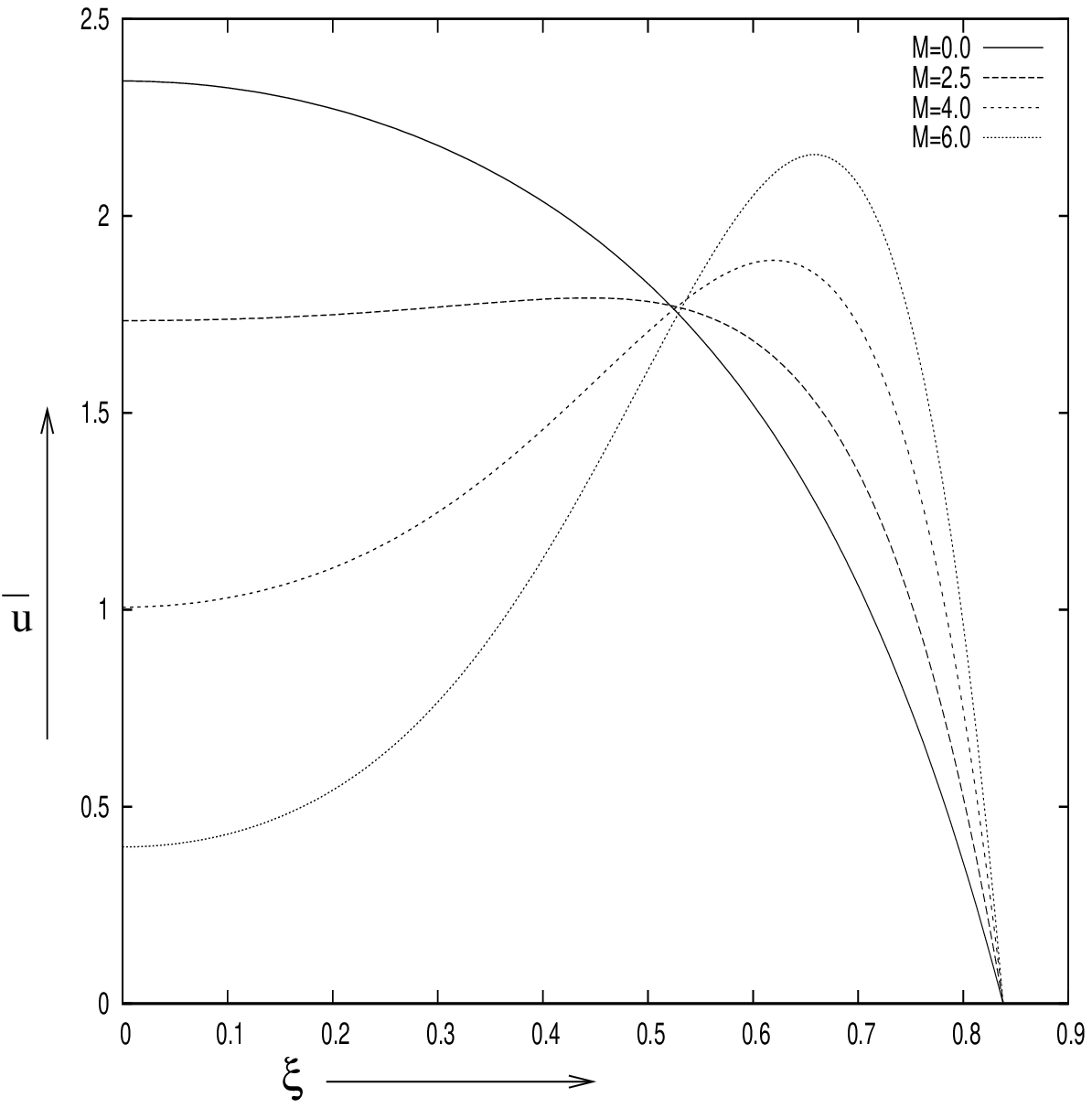}\\
      Fig. 6 Velocity distribution at $z=2.0$ for different values of $M$, when $H=0.2$, $\beta=2.5$, $\alpha=0.09$ and $k=0.25$.\\
\end{center}
\end{figure}
\begin{figure}[!ht]
\begin{center}
             \includegraphics[width=4.6in,height=3.5in ]{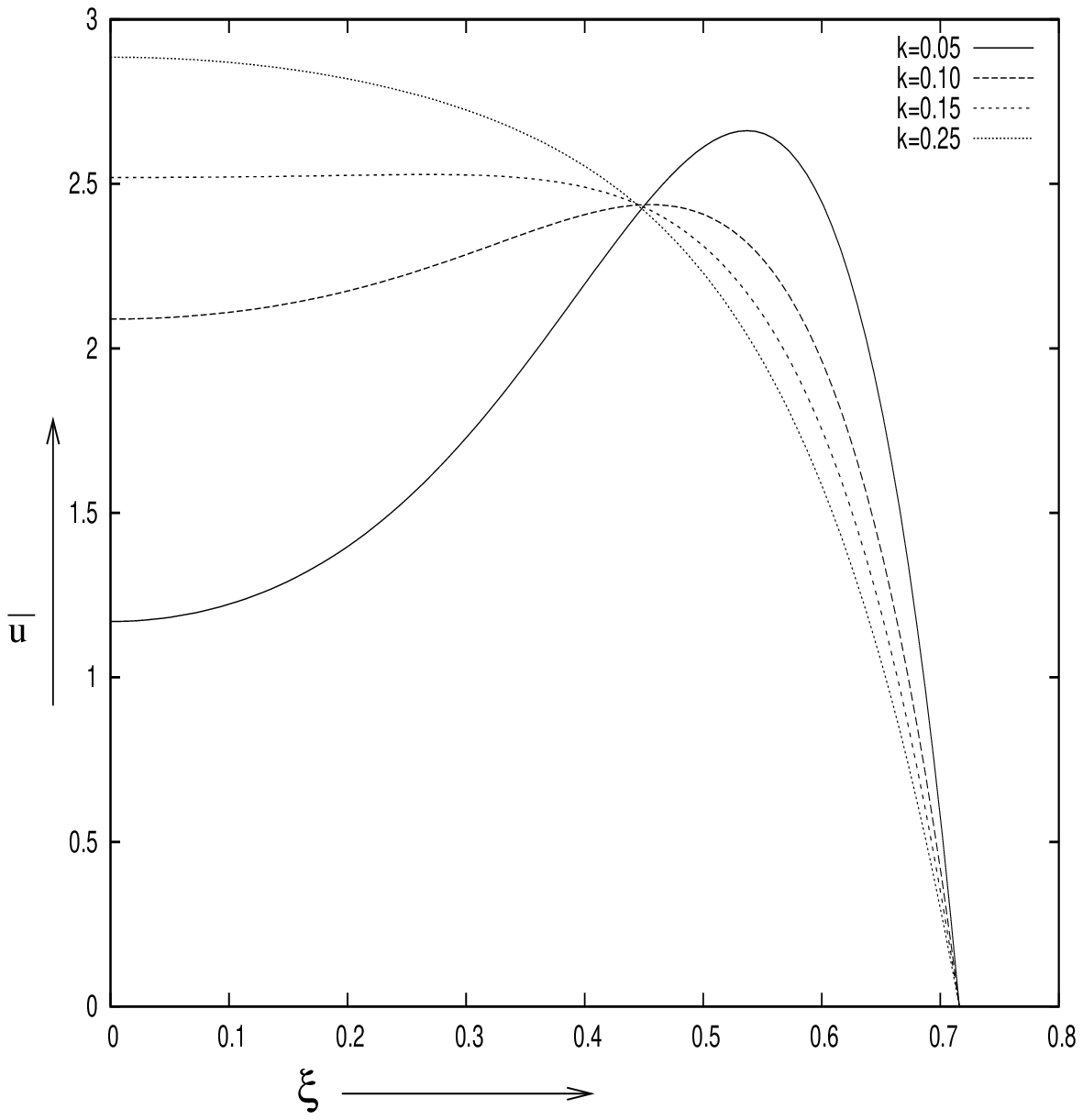}\\
      Fig. 7  Variation of axial velocity for different values of the permeability parameter $k$ at $z=2.0$, when $H=0.2$, $M=2.5$, $\beta=2.5$ and $\alpha=0.09$.\\

          \includegraphics[width=4.6in,height=3.5in ]{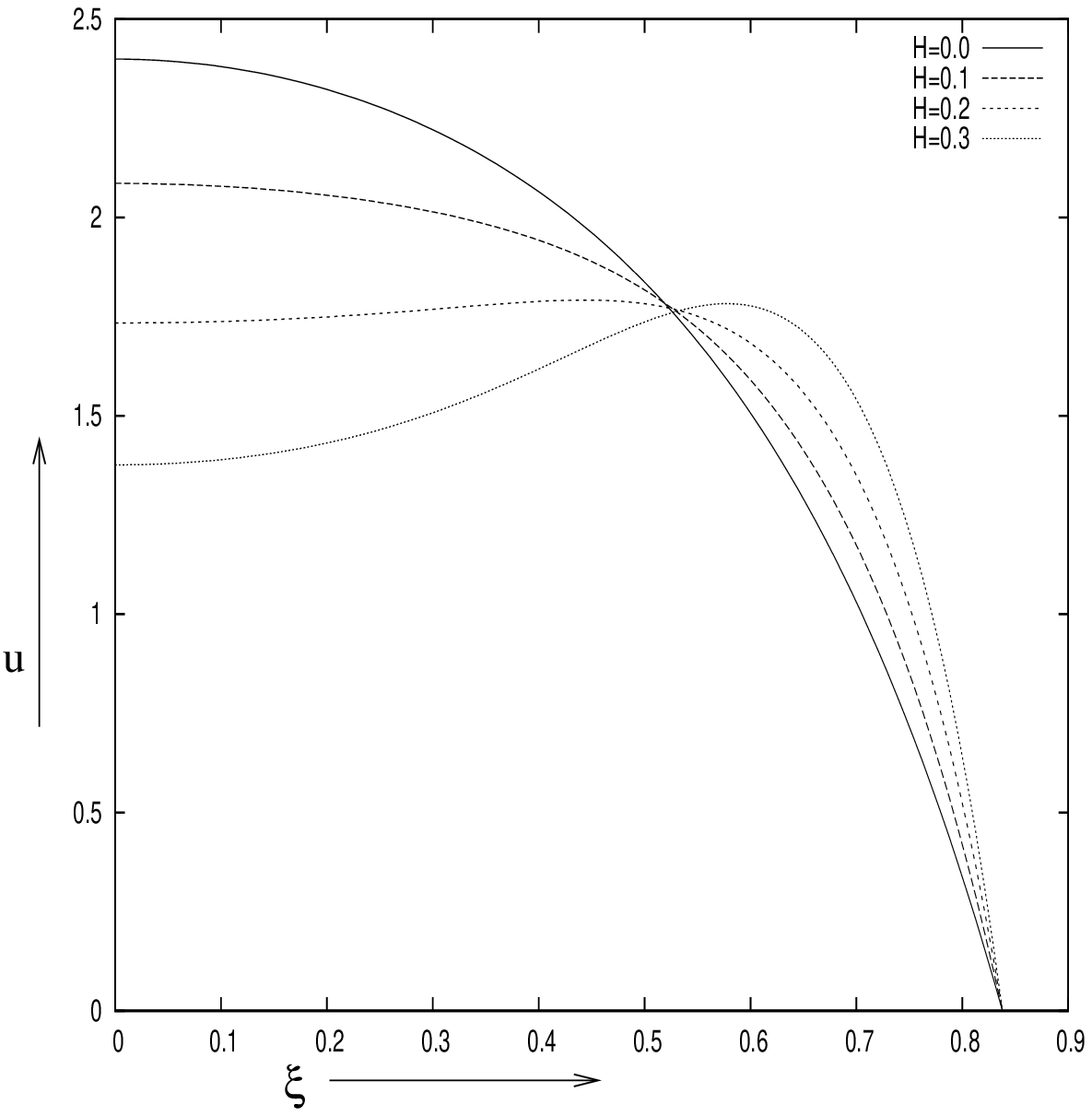}\\
       Fig.8 Variation of axial velocity in the radial direction for different values of $H$ $z=2.0$, when $M=2.5$, $\beta=2.5$, $\alpha=0.09$ and $k=0.25$.\\
\end{center}
\end{figure}
\begin{figure}[!ht]
\begin{center}
             \includegraphics[width=4.6in,height=3.5in ]{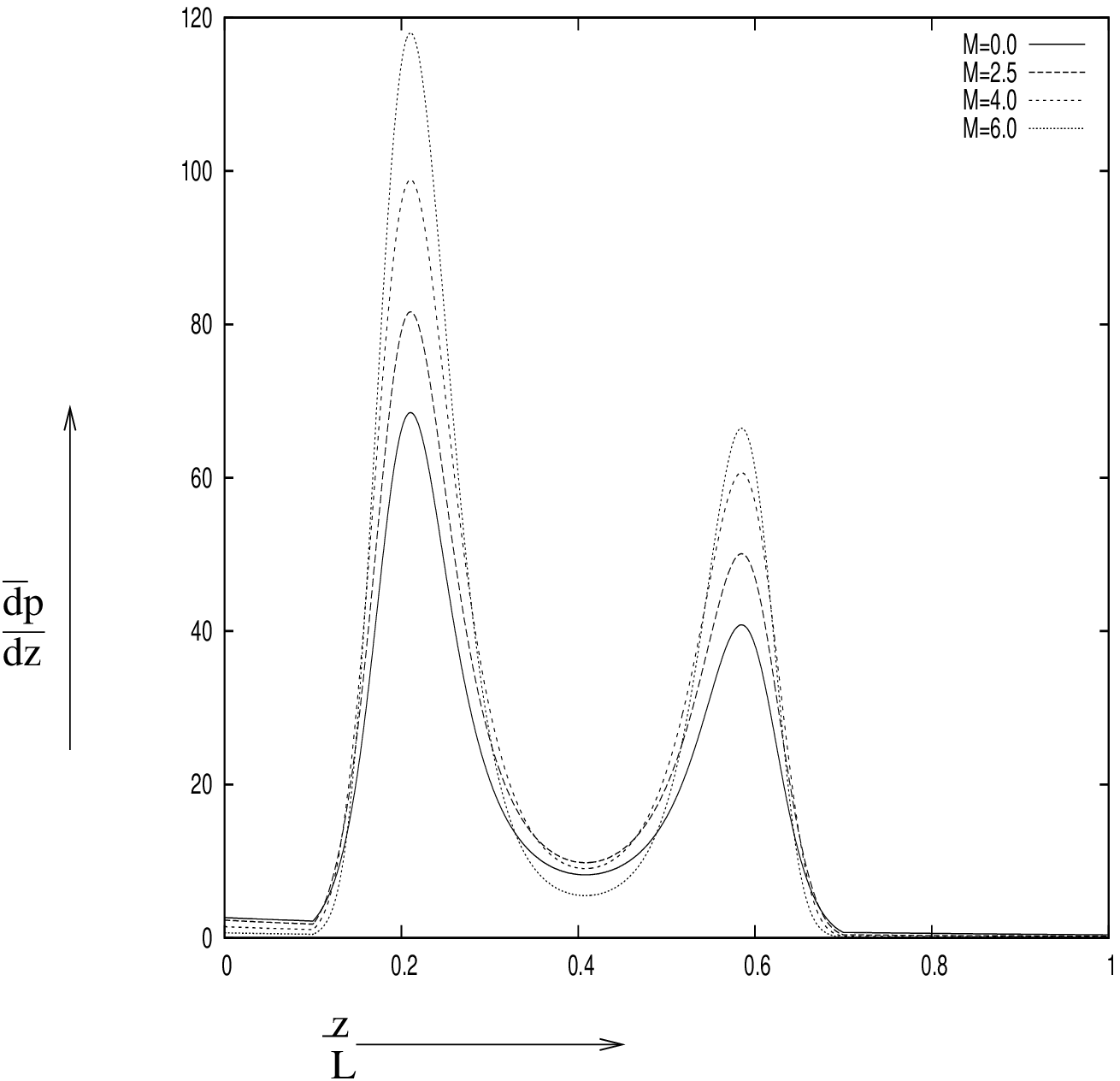}\\
      Fig. 9 Variation of pressure gradient ($\bar{\frac{dp}{dz}}$) for different values of $M$  when $H=0.2$, $\beta=2.5$, $\alpha=0.09$ and $k=0.25$.\\

             \includegraphics[width=4.6in,height=3.5in ]{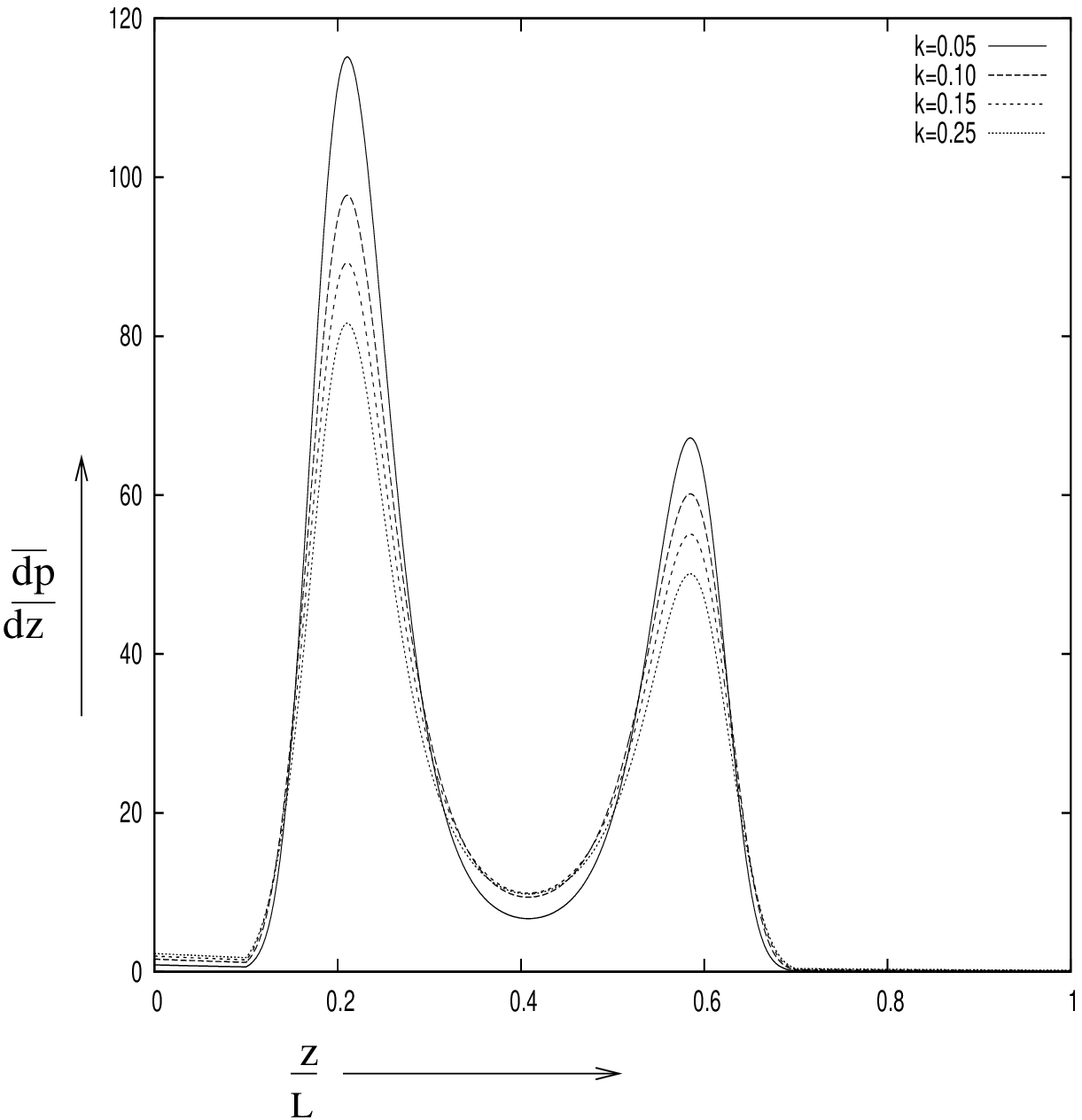}\\
      Fig. 10 Variation of pressure gradient ($\bar{\frac{dp}{dz}}$) with $z$ for different values of $k$ when $H=0.2$, $M=2.5$, $\beta=2.5$, when $\alpha=0.09$.\\
\end{center}
\end{figure}
\begin{figure}[!ht]
\begin{center}
             \includegraphics[width=4.6in,height=3.5in ]{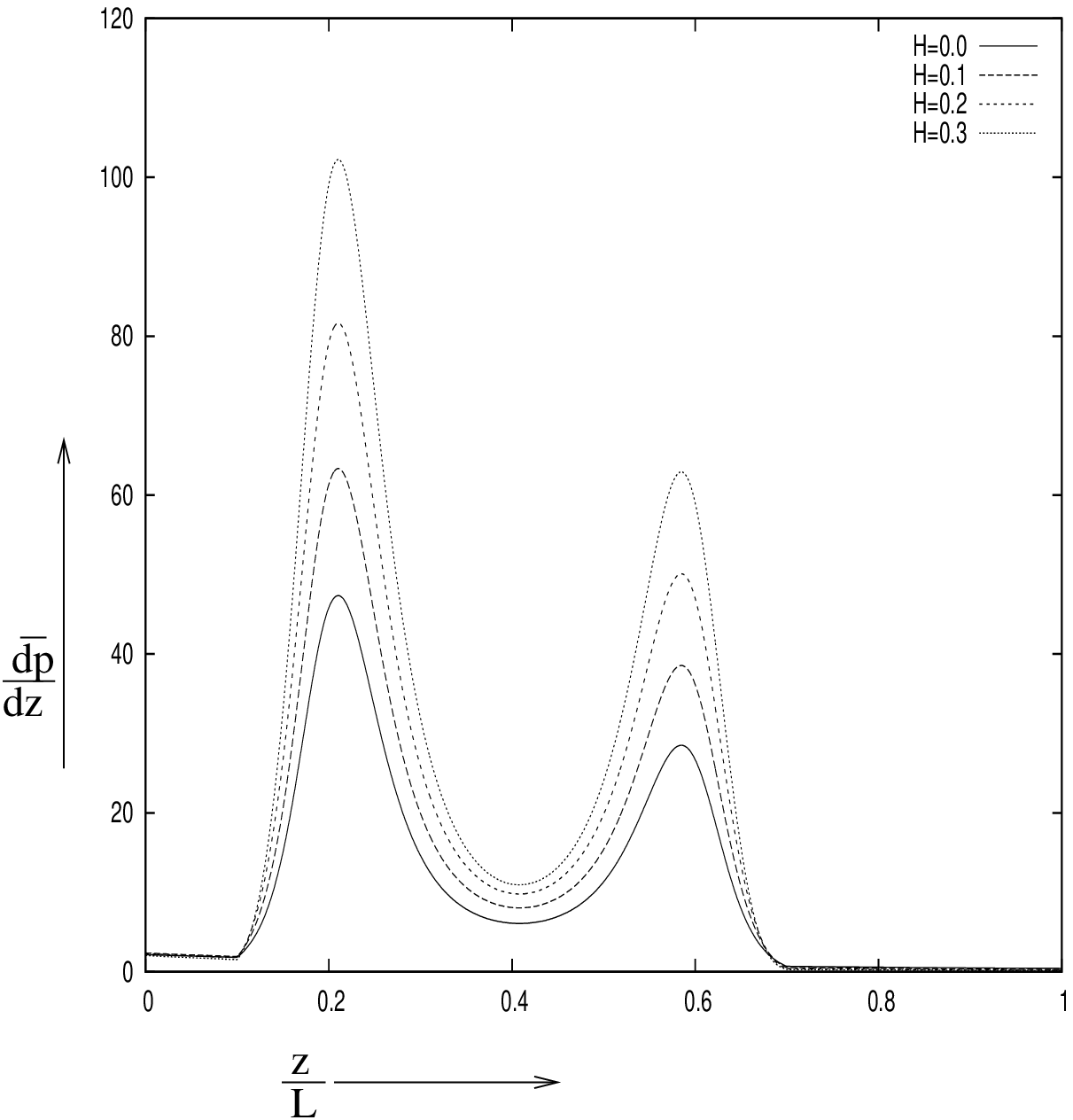}\\
      Fig. 11 Variation of pressure gradient ($\bar{\frac{dp}{dz}}$) for different values of $H$  when  $M=2.5$, $\beta=2.5$, $\alpha=0.09$ and $k=0.25$.\\

             \includegraphics[width=4.6in,height=3.5in ]{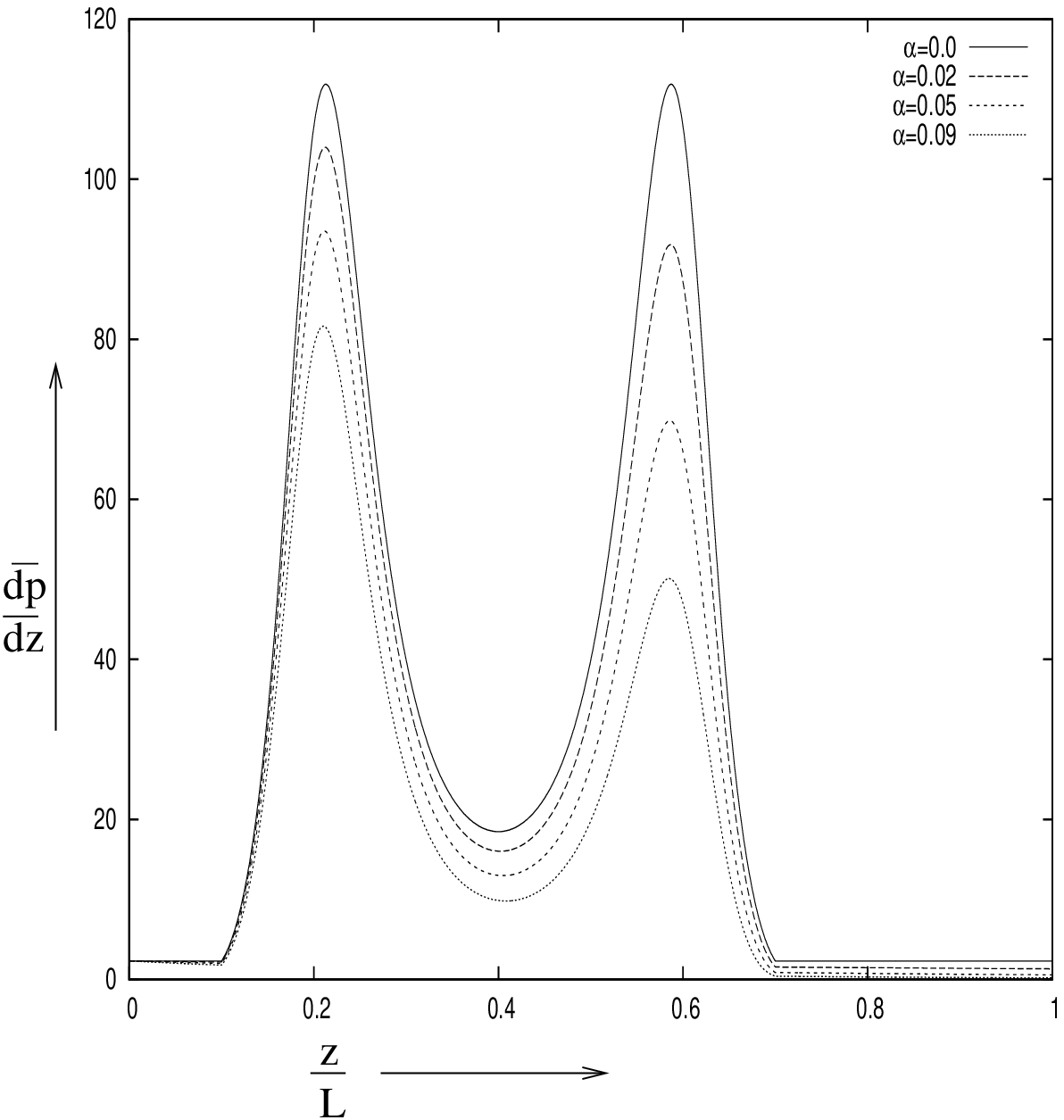}\\
      Fig. 12 Variation of pressure gradient ($\bar{\frac{dp}{dz}}$) with $z$ for different values of $\alpha$ when $H=0.2$, $M=2.5$, $\beta=2.5$, and $k=0.25$.\\
\end{center}
\end{figure}
\begin{figure}[!ht]
\begin{center}
          \includegraphics[width=4.6in,height=3.5in ]{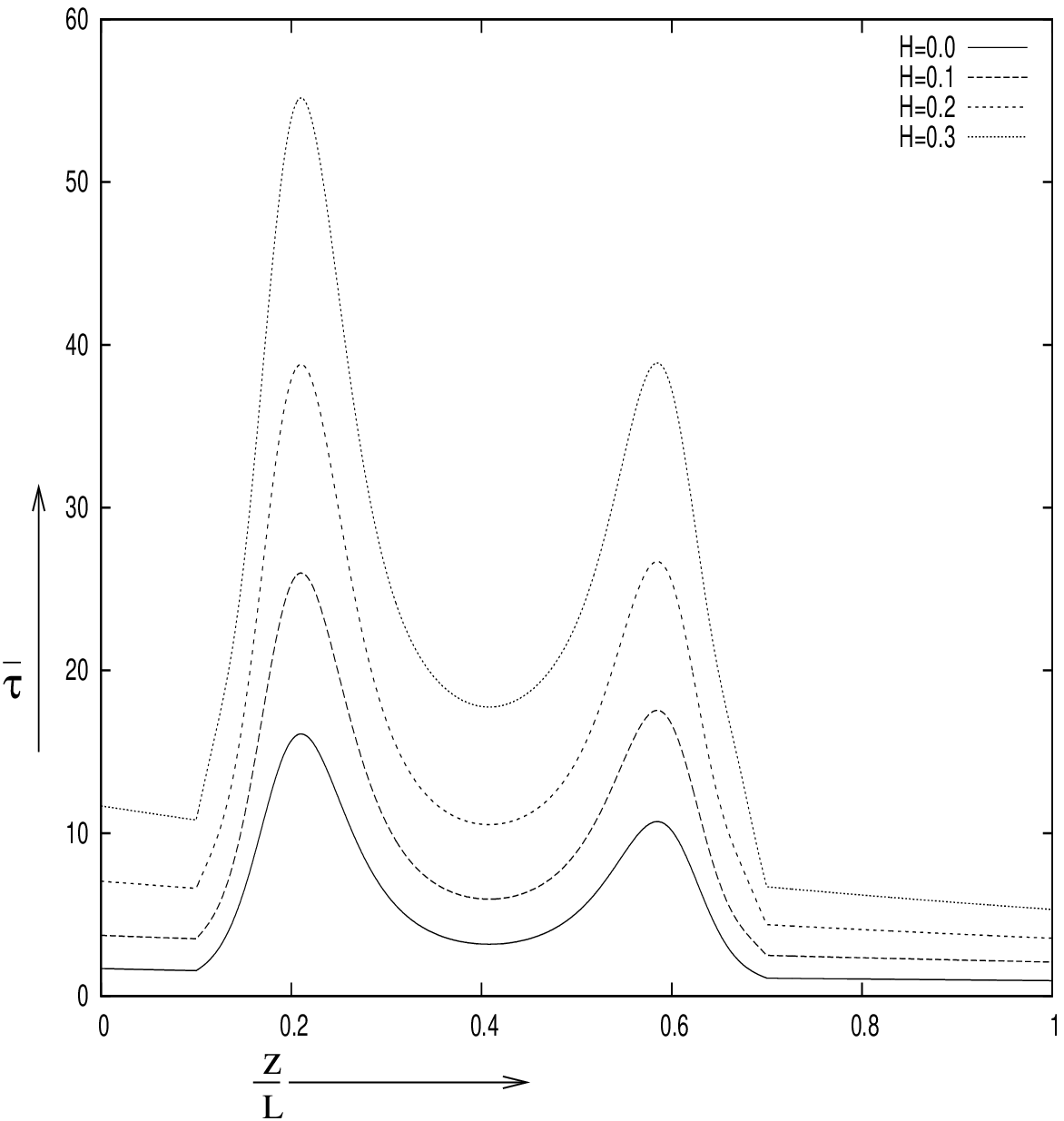}\\
      Fig. 13 Distribution of wall shear stress $\bar \tau$ along with $z$  for different values of $H$  when $k=0.25$, $M=2.5$, $\beta=2.5$, and $\alpha=0.09$\\
             \includegraphics[width=4.6in,height=3.5in ]{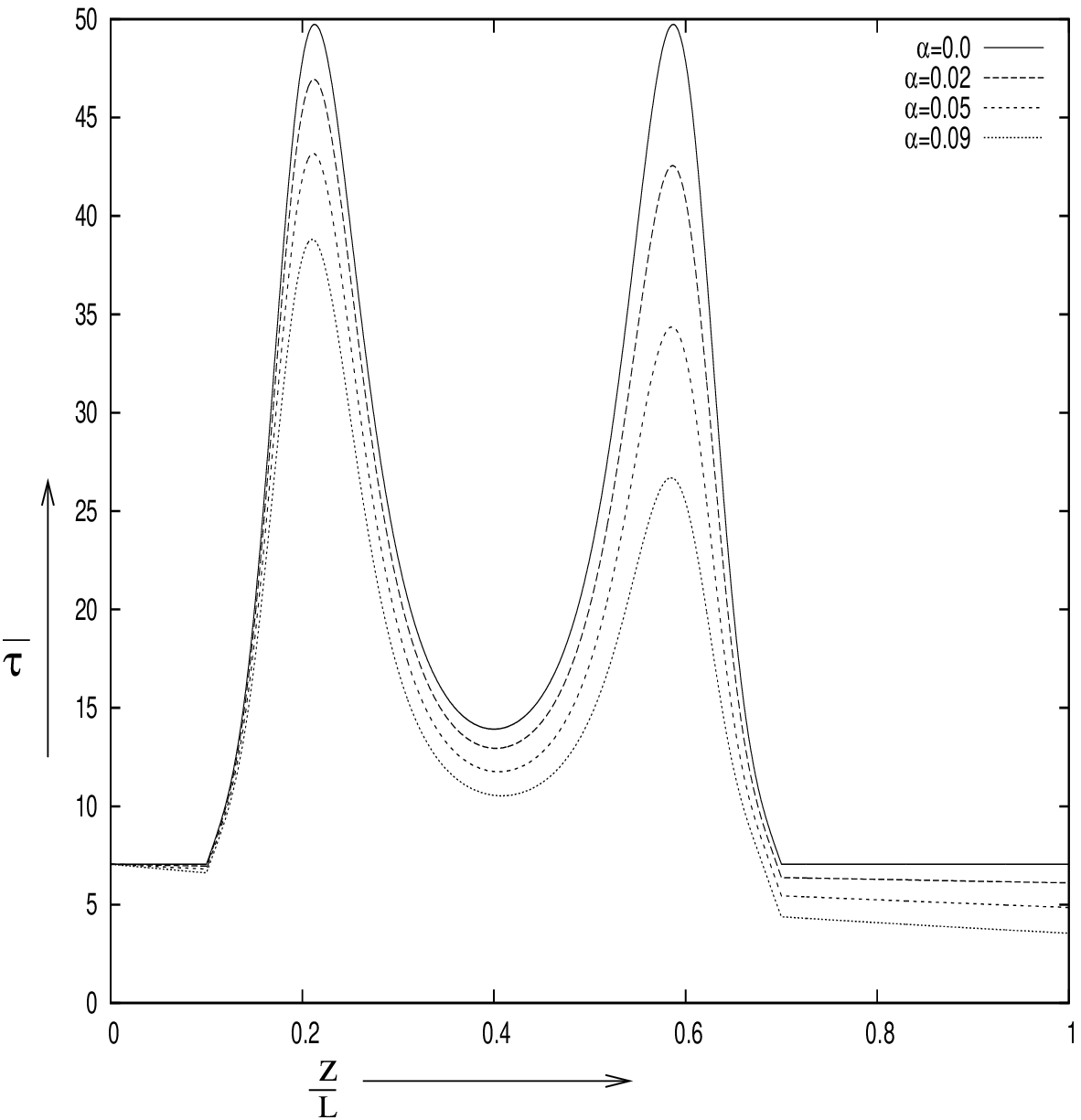}\\
      Fig. 14 Distribution of wall shear stress $\bar \tau$ along with $z$ for different values of $\alpha$  when $H=0.2$, $\beta=2.5$, $M=2.5$ and $k=0.25$.

\end{center}
\end{figure}
\end{document}